\magnification=\magstep1

\documentstyle{amsppt}
\NoBlackBoxes
\pagewidth{30pc}
\pageheight{50.5pc}

\def\Ad{\operatorname{ Ad}} 
\def\char{\operatorname{ char}}

\def\E{{\Cal E}}  
 
\def\F{{\Cal F}} 
\def\G{{\Cal G}} 
\def\har#1{\smash{\mathop{\hbox to .8 cm{\rightarrowfill}} 
\limits^{\scriptstyle#1}_{}}} 
 
\def\HNP{\operatorname{ HNP}}

\def\vert{\operatorname{{vert}}} 
\def\SL{\operatorname{SL} }

\def\L{{\Cal L}} 
\def\min{{\operatorname{ min}}} 
\def\max{{\operatorname{ max}}}

\def\R{{\Bbb R}} 
\def\GL{{\operatorname{ GL}}} 
\def\SL{{\operatorname{ SL}}}
 
\def\Int{\operatorname{ Int}}

\def\ot{\otimes}

\def\pr{\noindent{\sl Proof. }} 
\def\P{{\Cal P}} 
 
\def\JH{\operatorname{ JH}} 
\def\Q{{\Cal Q}}

\def\re{\noindent{\sl  Remark. }} 
\def\rk{\operatorname{ rk}}

\def\t{\tilde}

\def\PGL{\operatorname{ PGL}} 
\def\gl{\operatorname{\frak g\frak l} } 
\def\HN{\operatorname{ HN}} 
 
\def\Z{{\Bbb Z}}

\def\today{\ifcase\month \or January \or February \or March \or April 
\or May \or June \or July \or August \or September \or October 
\or November \or December \fi \space \number\day,  \number \year} 
\def\:{\colon} 
 
\topmatter 
\title  
Semistable principal $G$-bundles in positive characteristic  
\endtitle 
\author Adrian Langer 
\endauthor 
 
\rightheadtext{Semistable $G$-bundles} 
 

\address{ Adrian Langer:}\endaddress
\address{ 
1. Institute of Mathematics, Warsaw University,  Warszawa, Poland}\endaddress
\address{
2. Institute of Mathematics, Polish Academy of Science,
   Warszawa, Poland\medskip}
\endaddress

\email{alan\@mimuw.edu.pl}\endemail

\abstract  
Let $X$ be a normal projective variety defined over an algebraically
closed field $k$ of positive characteristic. Let $G$ be a connected
reductive group defined over $k$. We prove that some Frobenius pull
back of a principal $G$-bundle admits the canonical reduction $E_P$
such that its extension by $P\to P/R_u(P)$ is strongly semistable (see
Theorem 5.1).

Then we show that there is only a small difference between
semistability of a principal $G$-bundle and semistability of its
Frobenius pull back (see Theorem 6.3). This and the boundedness of the
family of semistable torsion free sheaves imply the boundedness of
semistable (rational) principal $G$-bundles.
\endabstract 
 
 
\endtopmatter 
 
\document  
\vskip .3 cm 

\heading 0. Introduction \endheading 
 
Let $k$ be an algebraically closed field of arbitrary
characteristic. Let $X$ be a normal projective variety over $k$ with a
very ample divisor $H$. One can then define the degree $\deg E$ of a
torsion free sheaf $E$ on $X$ with respect to $H$ and its slope $\mu
(E)=\deg E/\rk E$.  We say that a torsion free sheaf $E$ on $X$ is
{\sl slope semistable} with respect to $H$ if for every subsheaf $F$
of $E$ we have $\mu (F)\le \mu ( E).$ Every torsion free sheaf has a
canonical filtration with semistable quotients, the so called {\sl
Harder--Narasimhan filtration}. Let $\mu _{\max}(E)$ denote the slope
of the first factor of this filtration, i.e., the slope of the maximal
destabilizing subsheaf of $E.$

Assume that $\char k=p$ and let $F\: X\to X$ be the Frobenius
morphism.  It is well known that if $E$ is semistable then its
Frobenius pull back $F^*E$ need not longer be semistable. If all the
Frobenius pull backs $(F^k)^*E$ are semistable then $E$ is called {\sl
strongly semistable}.  By the Ramanan--Ramanathan theorem (see Theorem
2.7) such sheaves are well behaved under tensor operations. In
particular, a torsion free part of the tensor product of two strongly
semistable torsion free sheaves is strongly semistable.
 
In [La1] the author proved that for every torsion free sheaf $E$ there
exists some non-negative integer $l$ such that the factors of the
Harder--Narasimhan filtration of $(F^l)^*E$ are strongly
semistable. Therefore if $E_1$ and $E_2$ are torsion free sheaves
there exists some $l$ such that
$$\mu _{\max}((F^l)^*(E_1\ot E_2))=\mu _{\max}((F^l)^*E_1)+\mu
_{\max}((F^l)^*E_2).$$
On the other hand one can easily see that 
$$\mu _{\max}(E_1\ot E_2)\le {\mu _{\max}((F^l)^*(E_1\ot E_2))\over p^l}.$$
One can also show that the differences
$${\mu _{\max}((F^l)^*E_1)\over p^l}-\mu _{\max}(E_1)\quad \hbox{ and
} \quad {\mu _{\max}((F^l)^*E_2)\over p^l}-\mu _{\max}(E_2)$$ are
non-negative and bounded from the above by some explicit numbers
depending only on $(X, H)$, $p$ and the ranks of $E_1$ and $E_2$ (see
Theorem 1.3). Therefore we get a precise bound on the degree of
instability of the tensor product with respect to the degree of
instability of $E_1$ and $E_2$. This implies, e.g., that for large
primes $p$ the tensor product of semistable sheaves is semistable.

The main aim of this paper is to develop an analogue of the above
approach in the case of principal $G$-bundles, or more generally,
rational $G$-bundles.

The first step in the above approach is to prove that some Frobenius
pull back of a principal $G$-bundle $E$ admits the strongly semistable
canonical reduction (see Theorem 5.1). This is proved using Behrend's
combinatorial method allowing to compare degrees of two different
parabolic subschemes of a reductive group scheme (see [Be2]).

Then, for a semistable $G$-bundle $E$, we bound the degree of the
canonical reduction of any Frobenius pull back $(F^k)^*E$ (see
Corollary 6.6). The proof is similar to the proof of [La1], Corollary 6.2.
This result can be used to bound the degree of the associated bundle
$E(\frak g)$ of Lie algebras.

This can be used to reduce the problem of boundedness of the family
of semistable principal $G$-bundles with fixed numerical data to the
boundedness of torsion free sheaves. Then as a corollary of [La1],
Theorem 4.2, we get the following theorem (see Theorem 7.4):

\proclaim{Theorem 0.1} Let $G$ be a connected reductive group over $k$. 
Let us fix a polynomial $P$ and some constant $C$. Then the family of
all (rational) principal  $G$-bundles $E$ on $X$ such that the degree
of the canonical parabolic of $E$ is $\le C$, the degree of $E$ is
fixed, and the Hilbert polynomial of a torsion free sheaf extending
$E(\frak g)$ is equal to $P$, is bounded.
\endproclaim

In particular, the above theorem implies that the family of semistable
principal (rational) $G$-bundles with fixed degree and the Hilbert
polynomial is bounded. In the case when $X$ is a curve over an
algebraically closed field of characteristic zero, the boundedness of
the family of semistable principal $G$-bundles was proved by
A. Ramanathan (see [Ra2] and [Ra3]). His method works also in the
higher dimensional case.

In the case when $X$ is a curve over an algebraically closed field of
positive characteristic, the boundedness was first proved by
K. Behrend in his thesis (see [Be1], Theorem 8.2.6), using Harder's
result on reductions of $G$-bundles to a Borel subgroup of $G$. The
same proof was later published by Y. Holla and M. S. Narasimhan (see [HN]).

Once we have Corollary 6.6 our method is quite similar to Ramanathan's
method (with additional complications caused by higher dimension and
positive characteristic).

Now assume for simplicity that $G$ is semisimple and let $\rho \: G\to
\SL (V)$ be a homomorphism. Then the above results also allow us to bound 
the slope of the maximal destabilizing subsheaf of $E(V)$ for a
semistable $G$-bundle $E$ (see Theorem 8.4). In particular, if the
characteristic of the field is large and $E$ is semistable then we
prove that $E(V)$ is also semistable.

A similar theorem, but with usually better bounds on the
characteristic, was proved by S. Ilangovan, V. B. Mehta and
A. J. Parameswaran (see [IMP]). Our approach has the advantage of
giving some information for small primes.

In the forthcoming paper we will show how to apply the above results
and methods of [La1] to obtain restriction theorems for principal
$G$-bundles (see [BG] for some restriction theorems in the
characteristic zero case).

There is a recent preprint by F. Coiai and Y. Holla (see {\tt
math.AG/0312280}) in which the authors prove a non-effective version
of our Theorem 8.4 in the case when $H=\GL (V)$. They do it refining
the methods of [RR] and, to the author, it does not seem possible to
obtain effective bounds on instability of associated bundles using
their method.  They also try to use the above result to prove a weak
version of boundedness for semistable principal $G$-bundles defined
over the whole variety.

\medskip
The structure of the paper is as follows. In Section 1 we recall a few
results needed in the following. In Section 2 we define and study the
canonical reduction. In Section 3 we define the strong canonical
reduction and we study its properties.  In Section 4 we explain a
geometric meaning of complementary polyhedra introduced by Behrend in
[Be2]. In Section 5 we prove that some Frobenius pull back of a
principal $G$-bundle admits the strong canonical reduction.  In
Section 6 we study differences between semistability of a principal
$G$-bundle and of its Frobenius pull back. We apply these results in
Section 7 to get the boundedness of semistable principal
$G$-bundles. In Section 8 we show how to bound the degree of
instability of extensions of semistable $G$-bundles.

\bigskip

{\sl Notation.} 

We fix some notation used throughout the paper.  Let $k$ be an
algebraically closed field of arbitrary characteristic. Let $X$ be a
normal projective geometrically connected variety defined over $k$.
An open subset $U$ of $X$ is called {\sl big}, if the complement of
$U$ has codimension $\ge 2$.  A {\sl rational vector bundle} $E$ is a
vector bundle defined over some big open subset of $X$. In this case
$E$ (or, more precisely, the associated locally free sheaf) has a
unique extension $\t E$ to a reflexive sheaf on $X$.

Let $d$ be the dimension of $X$ and let $H_1\dots , H_{d-1}$ be ample
divisors on $X$. Then we can define the degree of $E$ as the degree of
its extension $\t E$ with respect to $H_1\dots H_{d-1}$. Using it one
can easily define slopes and stability of rational vector bundles
(cf.~[La2], Appendix).  In this paper, unless otherwise stated, we
will always talk about semistability defined with respect to the above
$1$-cycle.

Let $G$ be a connected reductive group over $k$. Then we define a {\sl
rational $G$-bundle} as a principal $G$-bundle on a big open
subset of the smooth part of $X$.  In this paper when writing ``a
principal $G$-bundle'' we will always mean only a rational $G$-bundle.

Let $R(G)$ denote the radical of $G$. Since $G$ is reductive, $R(G)$
is equal to the identity component of the reduced centre of $G$.

\heading 1. Preliminaries \endheading 
 
{\sl 1.1.} Here we recall some basic facts about parabolic subgroups
in reductive groups.

Let $G$ be a connected reductive algebraic group over $k$ and let
$\frak g$ denote its Lie algebra. Let us fix a maximal torus $T$ in
$G$ and a Borel subgroup $B$ containing it.  Let $X^*(T)$ be the
character group of $T$ and let $\Phi=\Phi (G,T)\subset X^*(T)$ be the
set of roots of $G$ with respect to $T$. By definition $\Phi$ is the
set of non-zero weights of $T$ in $\frak g$, acting via the adjoint
representation $\Ad$.  The choice of $B$ determines the set $\Phi ^+$ of
positive roots which contains the subset $\Delta$ consisting of simple
roots.
 
For any root $\alpha \in \Phi$ there exists an isomorphism $x_\alpha$
of $ G_a$ onto a unique closed subgroup $ X_\alpha $ of $ G$ such that
for any $t\in T$ and $a\in G_a$ we have $tx_\alpha (a)t^{-1}=x_\alpha
(\alpha (t)a)$.
 
There is a $1-1$ correspondence between subsets $I$ of $\Delta$ and
parabolic subgroups $P_I$ containing $B$. There are two possible
choices to define this correspondence. We do it in such a way that the
Levi subgroup $L_I$ of $P_I$ containing $T$ is generated by $T$ and
$X_{\pm \alpha}$ for $\alpha \in \Delta -I$. Then $B$ corresponds to
$\Delta$ and $G$ corresponds to $\emptyset$.

Let us fix $I$ and let $\Phi _I$ be the subset of $\Phi ^+$ consisting of those roots
that are linear combinations of roots in $\Delta -I$. Any root $\alpha
\in \Phi ^+ -\Phi _I$ can be written as  
$\alpha =\sum _{\alpha _i \in \Delta } n_{i} \alpha _i $, where $n_i
\ge 0$ and $l(\alpha)=\sum _{\alpha _i \in I}n_i >0$.  The number
$l(\alpha)$ is called the {\sl level} of $\alpha$ and $S(\alpha )=\sum
_{\alpha _i \in I}n_{i} \alpha _i$ is called the {\sl shape} of
$\alpha$.
 
For each non-zero shape $S$ we set $V_S=\prod _{S(\alpha )=S}
X_{\alpha}$.  Each $V_S$ is a module over the Levi subgroup $L=L_I$ of
$P=P_I$ acting by inner automorphisms. One can also see that $R(L)$
acts on $V_S$ by scalars.
 
The unipotent radical $R_u (P)$ is generated by $X_\alpha$ for $\alpha
\in \Phi ^+- \Phi _I$ and it has a natural filtration $U_m\subset
\dots \subset U_1\subset U_0=R_u(P)$ such that $U_i\triangleleft R_u(P)$.
$U_i$ is defined as $\prod _{l(\alpha)>i}X_\alpha$.  For each factor
of this filtration we have the decomposition $U_i/U_{i+1}=\bigoplus
_{l(S)=i+1} V_S$ into a direct sum of $L$-modules. If $G$ is not
special then each $V_S$ is a simple $L$-module and the filtration is
the socle (Loewy) series of $R_u (P_I)$, where $R_u (P)$ is treated as
a $P$-module with $P$ acting by inner automorphisms (see [ABS], Lemma
4). However, if $G$ is special it can happen that $V_S$ is not a
simple $L$-module (see [ABS], Section 3, Remark 1).
 
Let $\frak g _{\alpha}$ be the Lie algebra of $X_\alpha$. Then we have
an induced filtration of the Lie algebra $\frak u$ of $R_u(P)$ with quotients being
$L$-modules, and we can identify the corresponding $L$-module $\oplus
_{S(\alpha)=S}\frak g_\alpha$ with $V_S$.

Then one can see that $\frak g/\frak p$, where $\frak p$ is the Lie
algebra of $P$, has a dual filtration $W_m\subset \dots \subset
W_1\subset W_0=\frak g/\frak p$ such that $W_i/W_{i+1}=\bigoplus
_{l(S)=i+1} V_S^*=\bigoplus _{l(S)=-(i+1)} V_S$ (cf. [ABS], Section 3,
Remark 6).
 
\medskip
{\sl 1.2.} We also need to recall some notation from [La1]. If $E$ is
a rational vector bundle on $X$ defined over a field of characteristic
$p$ then we set
$$L_{\max}(E)=\lim _{k\to \infty} {\mu _{\max} ((F^k)^*E)\over p^k}.$$
This is a well defined rational number (this follows from [La1],
Theorem 2.7). Similarly, one can also define $L_{\min}(E)$.

Note that in the introduction we proved that
$$L_{\max}(E_1\ot E_2)=L_{\max} (E_1)+L_{\max}(E_2)$$
for any two rational vector bundles $E_1$ and $E_2$.

\proclaim{Theorem 1.3} {\rm ([La1], Corollary 6.2)}
Let $E$ be a rational vector bundle of rank $r$.  
\itemitem{(1)} If $\mu _{\max
}(\Omega _X)\le 0$ then $L_{\max}(E)=\mu _{\max} (E)$ and
$L_{\min}(E)=\mu _{\min} (E)$.  
\itemitem{(2)} If $\mu _{\max }(\Omega _X)> 0$ then
$$L_{\max}(E)\le \mu _{\max}(E)+{r-1\over p}L_{\max} (\Omega _X)$$
and
$$L _{\min}(E)\ge \mu _{\min}(E)-{r-1\over p}L_{\max} (\Omega _X).$$
\endproclaim

In Sections 6 and 8 we will prove similar theorems for principal
$G$-bundles.

\medskip

{\sl 1.4.} Let $X$ be a $d$-dimensional normal variety defined over an
algebraically closed field $k$ and let $H$ be an ample divisor on
$X$. Let $E$ be a rank $r$ torsion free sheaf on $X$. Then there exist
integers $a_0(E), \dots ,a_d (E)$ such that
$$\chi (X, E(mH))=\sum _{i=1}^d a_i (E){m+d-i \choose d-i}.$$ 

\proclaim{Theorem 1.5} {\rm ([La1], Theorem 4.4)}
Let $\mu _{\max}, a_0, a_1 $ and $a_2$ be some fixed numbers.  Then the
family of torsion free sheaves on $X$ such that $\mu _{\max}(E)\le \mu
_{\max}$, $a_0(E)=a_0$, $a_1(E)=a_1$ and $a_2(E)\ge a_2$ is bounded,
i.e., there exists a scheme $S$ of finite type over $k$ and an
$S$-flat sheaf $\E$ on $X\times S$ such that each member of the above
family is contained in $\{\E_s\}_{s\in S}$, where $\E _s$ is the
restriction of $\E$ to the fibre of the canonical projection over $s\in
S$.
\endproclaim

Since $a_0(E)=rH^d$ and $a_1(E)=(c_1E-{r\over 2}K_X)H^{d-1}$, fixing
$a_0(E)$ and $a_1(E)$ is equivalent to fixing the rank $r$ of $E$ and
the degree $c_1(E)H^{n-1}$ of $E$.

In the case $X$ is smooth we can define the discriminant $\Delta (E)$ of
$E$ as $2rc_2-(r-1)c_1^2$.  It is easy to see that the condition
$a_2(E)\ge a_2$ is equivalent to the condition $\Delta (E)H^{d-2}\le C_X (r,
a_1, a_2)$ for some explicit function $C_X$ depending only on $X$ and
$H$ (it is also equivalent to bounding $c_2(E)H^{d-2}$ from the above).

\heading 2. Harder--Narasimhan filtration \endheading

{\sl 2.1.} Let $G$ be a connected reductive group over $k$ and let $E$ be a
rational $G$-bundle on $X$.  Let $E(G)=E\times _{G, \Int}G$ denote the
group scheme associated to $E$ by the action of $G$ on itself by inner
automorphisms. Then we define the degree of $E(G)$ as the degree of
the Lie algebra bundle $E(\frak g)=E\times _{G, \Ad }\frak g$ of
$E(G)$ on $X$ considered as a rational vector bundle on $X$.
 
Let us fix a maximal torus $T$ in $G$ and some Borel subgroup $B$
containing $T$.  Let $P$ be a parabolic subgroup of $G$ containing $B$
and let $E_P$ be a (rational) reduction of its structure group to a
parabolic subgroup $P$. Since every parabolic subgroup of $G$ is
conjugate to exactly one parabolic subgroup of $G$ containing $B$, we
do not restrict the class of considered reductions (cf. [Ra2], Remark
3.5.7).
 
Let $E_P(P)=E_P \times _{P, \Int} P$ be the parabolic subgroup scheme 
of $E(G)$. 
 
\proclaim{Definition 2.2}  

The reduction $E_P$ of $E$ is called {\rm canonical} (or the 
{\rm Harder--Narasimhan filtration}) if it satisfies the following 
conditions: 
\item{(1)} for any parabolic subgroup scheme $Q\subset E(G)$ we have 
$\deg Q\le \deg E_P(P)$, 
\item{(2)} $E_P(P)$ is maximal among all parabolic subgroup schemes $\P$ of $E(G)$ that
satisfy (1),i.e., if $\P$ satisfies (1) and contains $E_P(P)$, then $\P =E_P(P)$. 
\endproclaim

The degree $\deg E_P (P)$ of the canonical reduction is denoted by
$\deg _{\HN}E$. This is a well defined integer.  This follows from the
fact that if $E_P$ is a reduction of the structure group of $E$ to $P$
then $\deg E_P (P)\le \rk E(\frak g) \cdot \mu _{\max} (E(\frak g))$,
since $E_P(\frak p)\subset E(\frak g)$ (this proof works in general;
cf.~[Be2], Lemma 4.3 for the curve case). Note that $\deg _{\HN}E\ge
0$, since $G$ is also parabolic and $\deg E(G)=0$.

\proclaim{Definition 2.3}

$E$ is called {\rm slope semistable} if and only if $\deg
_{\HN}E=0$, i.e., if the degree of any parabolic subgroup scheme of $E(G)$
is non-positive. 

$E$ is called {\rm strongly slope semistable} if and only if $\char
k=0$ or $\char k>0$ and $(F^l)^*E$ is slope semistable for all $l\ge
0$, where $F$ denotes the Frobenius morphism.
\endproclaim

The above definitions are not completely standard if one tries to
understand them in the vector bundle case. In this case one can easily
interpret the definition in the following way.  To any sheaf $G$ we can
associate the point $p(G)=(\rk G, \deg G)$ in the plane.  Let
$0=E_0\subset E_1\subset \dots \subset E_m=E$ be a filtration of a
vector bundle $E$ with torsion free quotients. We can successively
connect the points $p(E_0),\dots , p(E_m)$ by line segments. If we
also connect $p(E_m)$ with $p(E_0)$ then we obtain a generalized
polygon. One can easily see that the area of this polygon is equal to
half of the degree of the corresponding parabolic subscheme of the
associated group scheme.

The polygon corresponding to the Harder--Narasimhan filtration is
called the {\sl Harder--Narasimhan polygon} and it is denoted by $\HNP
(E)$.  Now Definition 2.2 says that the Harder--Narasimhan filtration
corresponds to the polygon with the largest area. This is clear, since
the Harder--Narasimhan polygon lies over all the polygons obtained
from the filtrations of $E$.  This and a small computation imply the
following proposition:

\proclaim{Proposition 2.4}
Let $E$ be a rational $\GL (V)$-bundle and let $E(V)$ be the
corresponding rational vector bundle. Let $r_i, \mu _i$ denote ranks
and slopes of the quotients of the Harder--Narasimhan filtration of
$E(V)$.  Let $\mu _{\max}, \mu _{\min}$ denote the corresponding
slopes for $E(V)$.  Then
$$\deg _{\HN} E=2\operatorname{area}\, \HNP (E(V)) =\sum _{i<j} r_ir_j(\mu _i-\mu _j).$$
In particular, we have
$$(r-1)(\mu _{\max}-\mu _{\min })\le \deg _{\HN} (E)\le {r^2\over 4}(\mu _{\max}-\mu _{\min }),$$
where $r=\dim V$. 
\endproclaim

\medskip 
If $F$ is a vector bundle then $\deg _{HN}F$ will denote
$2\operatorname{area}\, \HNP (F)$, which by the above proposition is
equal to the degree of the canonical parabolic of the corresponding
principal bundle.

\medskip
If $E$ is a principal $G$-bundle and $k$ is a field of characteristic $p$ then we set   
$$\deg _{\HN , l}E={\deg _{\HN} (F^l)^*E\over p^l}.$$  
 
\proclaim{Lemma 2.5} 
The sequence $\{\deg _{\HN , l}E\}$ is non-decreasing and it 
has a finite limit denoted by $\deg _{\HN , \infty}E$.  
\endproclaim 
 
\pr The first part follows immediately from the definition of the canonical  
reduction.  To prove existence of the limit it is sufficient to note
that if ${\tilde E}_P$ is the canonical reduction of $(F^k)^*E$ then
the vector bundle ${\tilde E}_P(\frak p)$ is contained in the vector
bundle $(F^k)^*E(\frak g)$, so its degree is less or equal to $p^k \rk
E(\frak g)\cdot L_{\max} (E(\frak g))$ (see [La], Corollary 2.5). In
particular, $\deg _{\HN , \infty}E\le \rk E(\frak g)\cdot L_{\max}
(E(\frak g)),$ Q.E.D.
 
\medskip

Let $L$ be a Levi subgroup of $P$. We have a natural projection $P\to 
P/R_u(P)$, where $R_u(P)$ is the unipotent radical of $P$. By the 
definition of a Levi subgroup we have $P/R_u(P)\simeq L$, so each principal 
$P$-bundle has an extension to an $L$-bundle. 
 
Let us also recall that the unipotent radical $R_u(P)$ has the
filtration $U_m\subset U_{m-1}\subset \dots \subset U_1=R_u(P)$ in
which the quotients $U_i/U_{i+1}$ are direct sums of simple
$L$-modules $V_S$ for all shapes of level $i+1$ (see 1.1). Each
$V_S$ is also a $k$-vector space.

\proclaim{Theorem 2.6} {\rm ([Be2], Theorem 7.3)} 
Every principal $G$-bundle has a canonical reduction $E_P$ to some 
parabolic subgroup $P$ of $G$.  This reduction satisfies the following 
conditions: 
\item{(1)} the extension $E_L$ of $E_P$ to $L$ is a semistable rational  
$L$-bundle,
\item{(2)} for all shapes $S$ of positive level the associated vector 
bundle $E_L(V_S)$ has positive degree.

Moreover, if (1) and (2) are satisfied for some reduction $E_P$ of $E$ 
to some parabolic subgroup $P$ of $G$, then this reduction is 
canonical. 
\endproclaim 
 
The above theorem is formulated in [Be2] only for principal $G$-bundles
(or, more generally, group schemes) over a curve. However, its proof
with minor modifications works in general. One should only note that
$\deg _{HN}E$ is well defined and then repeatedly use the fact that
any generic section of a parabolic subscheme of $E(G)$ extends to some
big open subset of $X$ (cf.~[RR], Section 4). 

In [Be2] condition (2) is asserted only for shapes of level $1$, but it
is true for all shapes of positive level (see Lemma 4.5).

In the vector bundle case (1) of Theorem 2.6 corresponds to convexity
of the Harder--Narasimhan polygon and (2) corresponds to semistability
of quotients in the Harder--Narasimhan filtration.

\medskip
The following theorem was proved by Ramanan and Ramanathan:

\proclaim{Theorem 2.7} {\rm ([RR], Theorem 3.23; see also [La2], Theorem A.3)} 
Let $\rho \:G \to H$ be a homomorphism of connected reductive
$k$-groups and assume that $\rho (R(G))\subset R(H)$. Let $E_G$ be a
rational $G$-bundle and let $E_H$ be the rational $H$-bundle obtained
from $E$ by extension.  If $E_G$ is strongly semistable then $E_H$ is
also strongly semistable.
\endproclaim

\proclaim{Corollary 2.8}
Let $E$ be a principal $G$-bundle. If $E(\frak g)$ is semistable as a
vector bundle then $E$ is semistable. In particular, $E$ is strongly
semistable if and only if $E(\frak g)$ is strongly semistable.
\endproclaim

\pr
If $E$ is not semistable then $E(\frak g)$ is a degree zero vector
bundle which contains the vector bundle $E_P(\frak p)$ of degree
$>0$. Hence $E(\frak g)$ is not semistable.

If $E$ is strongly semistable then $\Ad\: G\to \GL (\frak g)$ maps the
radical of $G$ to the identity, so by Theorem 2.7 $E_{\GL (\frak g)}$ is
a strongly semistable principal $\GL (\frak g)$-bundle.  This implies
that the associated vector bundle is strongly semistable, Q.E.D.

\medskip

Let us recall that canonical reduction is functorial under separable
base change:

\proclaim{Lemma 2.9} {\rm (see [Be2], Corollary 7.4)} 
Let $\pi \: Y\to X$ be a finite separable morphism of normal 
projective varieties over $k$.  Let $E_P$ be the canonical reduction 
of a rational $G$-bundle $E$ defined over $X$. Then $\pi ^*E_P$ is the 
canonical reduction of $\pi ^*E.$ 
\endproclaim

The following proposition is an analogue of the Ramanan--Ramanathan
theorem but we do not need to assume strong semistability.  It immediately
implies Theorem 2.7 in the case $\rho$ is surjective.

\proclaim{Proposition 2.10}
Let $\rho\: G\to H$ be a surjective homomorphism of connected
reductive groups. Assume that the kernel group scheme of $\rho$ is
contained in the centre group scheme $Z(G)$ of $G$. Let $E$ be a
principal $G$-bundle with the Harder--Narasimhan filtration
$E_P$. Then the extension of structure group of $E_P$ to the image $Q$
of $P$ is the Harder--Narasimhan filtration of the extension $E_H$ of
structure group of $E$ to $H$.
\endproclaim

\pr
Since $\ker \rho \subset Z(G)\subset P$, we have $G/P \simeq (G/\ker
\rho)/ (P/\ker \rho)=H/Q$. Since $\frak g/ \frak p$ is the tangent space at 
$e$ to $G/P$, we have $\frak g/ \frak p\simeq \frak h/ \frak q$ and hence
$$\deg E_Q(\frak q)= - \deg E_Q(\frak h/\frak q)= - \deg E_P (\frak
g/\frak p)=\deg _{\HN} E.$$ If $E_{Q'}$ is a reduction of $E_H$ to a
parabolic subgroup $Q'\subset H$ and $P'=\rho ^{-1}(Q')$, then an
isomorphism $G/P'\simeq H/Q'$ induces a reduction $E_{P'}$ of $E$ to
the parabolic $P'$ (let us recall that such a reduction can be treated
as a rational section of $E(G/P')\to X$). By a similar
computation as above we have
$$\deg E_{Q'}(\frak q ')=\deg E_{P'}(\frak p ')\le \deg _{\HN}E= \deg E_{Q}(\frak q).$$
This shows that $E_Q$ satisfies condition (1) of Definition 2.2. Similarly one can check
that $E_Q$ satisfies condition (2), Q.E.D.

\medskip
As a corollary we see that in arbitrary characteristic a principal
$G$-bundle is semistable if and only if its extension to the adjoint
form $\Ad G$ is semistable (cf. Corollary 2.8). In the following we do
not use this fact.

\heading 3. Strong Harder--Narasimhan filtration\endheading

\proclaim{Definition 3.1}  
Let $E$ be a principal $G$-bundle. We say that $E_P$ is the {\rm
strong Harder--Narasimhan filtration} of $E$ (or that $E_P$ is the
{\rm strong canonical reduction} of $E$) if it satisfies condition (2)
of Theorem 2.6 and the extension $E_L$ is strongly semistable.
\endproclaim 
 
Obviously, if $k$ has positive characteristic then not every principal
$G$-bundle has a strong Harder--Narasimhan filtration. However, on
some special manifolds (e.g., with a globally generated tangent
bundle) every semistable $G$-bundle is strongly semistable and then
every principal $G$-bundle has a strong Harder--Narasimhan filtration
(see Corollary 6.4).

\bigskip

In [Be2], Conjecture 7.6, Behrend conjectures that if $E_P$ is the
canonical reduction of $E$ then 
$$h^0 (X, E_P (\frak g /\frak p))=0.$$ This conjecture was proved by
V. B. Mehta and S. Subramanian under the assumption that the
characteristic of the base field is large (see [MS], Corollary 3.6).

Here we prove that this result holds in an arbitrary characteristic
for some Frobenius pull back of $E$ (see Theorem 5.1):

\proclaim{Proposition 3.2}  
Assume that a principal $G$-bundle $E$ has the strong
Harder--Na\-ra\-sim\-han filtration $E_P$. Then
$$h^0 (X, E_P (\frak g /\frak p))=0.$$
\endproclaim

\pr 
Let us note that $E_P (\frak g /\frak p)$ has a filtration by vector
bundles in which the quotients are duals of vector bundles $E_L (V_S)$
for shapes of level $\ge 1$. Since the radical $R(L)$ acts on $V_S$ by
scalars and $E_L$ is strongly semistable, it follows that $E_L(V_S)$
is also strongly semistable (see Theorem 2.7). Since each $E_L(V_S)$
has positive degree, we have $\mu _{\max} (E_P(\frak g /\frak p))<0$
and in particular $E_P(\frak g /\frak p)$ has no sections, Q.E.D.

\medskip

\proclaim{Proposition 3.3} 
Assume that $E$ has strong Harder--Narasimhan filtration $E_P$. Let 
$$E_{-r}\subset \dots \subset E_{-1}\subset E_0\subset \dots \subset
E_s=E(\frak g)$$ be the Harder--Narasimhan filtration of $E(\frak g)$
indexed in such a way that $\mu (E_i/E_{i-1})<0$ for $i\ge 1$ and $\mu
(E_i/E_{i-1})\ge 0$ for $i\le 0$.  Then $E_0=E_P(\frak p)$ and
$E_{-1}=E_P(\frak u)$. In particular, $E_0/E_{-1}=E_L (\frak l)$ is
strongly semistable of degree $0$, $r=s$ and $E(\frak g)/E_0$ is isomorphic to
$E_{-1}^*$ as a rational vector bundle.
\endproclaim

\pr
Let us recall that $\frak u=\operatorname{Lie} R_u(P)$ is filtered with
$L$-modules $V_S$, where $S$ are shapes of level $\ge 1$. Then $\frak
g/\frak p$ is filtered with dual $L$-modules $V_S^*$. Hence $E(\frak
g/\frak p)=E(\frak g)/E_P(\frak p)$ has a filtration, whose quotients
are strongly semistable vector bundles $E_L(V_S)^*$ of negative
degree. In particular, $\mu _{\max} (E(\frak g)/E_P(\frak
p))<0$. Since $\mu _{\min} (E_0)\ge 0$, this implies that the map
$E_0\to E(\frak g)/E_P(\frak p)$ is zero, i.e., $E_0\subset E_P(\frak
p)$.

Note that $E_P(\frak p)/E_P(\frak u)=E_L (\frak l)$ is strongly
semistable and has degree $0$ (as a bundle of reductive Lie
algebras). Moreover, $E_P(\frak u)$ has a filtration with quotients
that are strongly semistable sheaves $E_P(V_S)$ of positive degree.
Hence $\mu _{\min} (E_P(\frak p))=0< \mu _{\max }(E(\frak g)/E_0)$,
which implies that the map $E_P (\frak p)\to E(\frak g)/E_0$ is
zero, i.e., $E_P(\frak p)=E_0$.

Now the map $E_{-1}\to E(\frak g)/E_P(\frak u)$ is zero since $\mu
 _{\min} (E_{-1})=\mu (E_{-1}/E_{-2})>\mu (E_0/E_{-1})\ge 0= \mu
 _{\max }(E(\frak g)/E(\frak u))$.  Note that $\mu _{\max }(E(\frak
 g)/E_{-1})=\mu (E_0/E_{-1})\le 0$. This follows from the fact that
 $E_0=E_P(\frak p)$ has a filtration with strongly semistable
 quotients, whose slopes are non-positive.  Since $\mu _{\min}
 (E_P(\frak u))>0$ it follows that the map $E_P(\frak u) \to E(\frak
 g)/E_{-1}$ is zero. Hence $E_{-1}=E_P(\frak u)$, Q.E.D.

\medskip

The above proposition should be compared with the construction in [AB].

\medskip

\proclaim{Corollary 3.4}
Let $E$ be a principal $G$-bundle which has strong
Harder--Narasimhan filtration $E_P$. Let $E(\frak g)$ be the vector
bundle associated to $E$ by the adjoint representation $\Ad \:G\to \GL
(\frak g)$. Then
$$(\dim \frak g+\dim \frak l)\cdot \deg _{\HN}E\le \deg _{\HN}E(\frak g)
\le 2\dim \frak g\cdot \deg _{\HN} E,$$
where $\frak l$ is the Lie algebra of a Levi subgroup of $P$.  
\endproclaim

\pr 
Draw the Harder--Narasimhan polygon $A=\HNP (E(\frak g))$. By
Proposition 3.3, $A$ is contained in the rectangle whose two vertices
are equal to $0$ and $p(E(\frak g))$ and one side contains points
$p(E_{-1})$ and $p(E_0)$.  Thus Proposition 2.4 implies that
$$ \deg _{\HN}E(\frak g)=2\operatorname{area}\, (A)\le 2\dim \frak
g\cdot \deg _{\HN} E.$$ Another inequality follows from the fact that
$A$ contains the trapezium with vertices $0, p(E(\frak g)), p(E_{0})$
and $p(E_{-1})$, Q.E.D.

\medskip

\proclaim{Corollary 3.5}  
Let $G$ and $H$ be connected reductive groups and let $\varphi \:G \to
H$ be a homomorphism. Let $E_G$ be a principal $G$-bundle and $E_H$ be
the principal $H$-bundle obtained from $E$ by extension.  Let $\phi\:
E_G(\frak g)\to E_H (\frak h) $ be the induced homomorphism of Lie
algebra bundles.  Assume that both $E_G$ and $E_H$ have strong
Harder--Narasimhan filtrations $E_P$ and $E_Q$, respectively. Then
$\phi (E_P(\frak p))\subset E_Q (\frak q).$
\endproclaim 
 
\pr 
By Proposition 3.3 we have $\mu _{\min} (E_P(\frak p))\ge 0>\mu _{\max}
(E_H (\frak h)/E_Q(\frak q)).$ Hence the map $E_P(\frak p)\to
E_H(\frak h)/ E_Q (\frak q)$ is zero, Q.E.D.

\medskip
\re
If $k$ is a field of characteristic zero then Atiyah and Bott showed
that the canonical reduction is functorial with respect to $\varphi$
(see [AB], Proposition 10.4). One cannot hope that this is true in the
positive characteristic case. However, from the result of Ilangovan,
Mehta and Parameswaran (see [IMP]) one can see that the canonical
reduction is functorial if the characteristic of the field is large
enough (weaker bounds can be obtained from Theorem 8.4). One need only
to take such $p=\char k$ that the kernel and the cokernel of $\frak
g\to \frak h$ are, as $G$-modules, the direct summands of $\frak g $
and $\frak h$,  and use Proposition 3.3 (or better [MS], Proposition 2.2).

\heading 4. Complementary polyhedra and elementary vector bundles\endheading

In this section we introduce complementary polyhedra and we prove some
auxiliary results about elementary vector bundles. The best place to
find the necessary definitions and basic properties is Behrend's paper
[Be2] or his PhD thesis [Be1].

\medskip

{\sl 4.1.} Let the notation be as in 1.1 and 2.1. Let $V$ be a real vector
subspace of $X^*(T)\ot _{\Z} \R $ spanned by the set $\Phi=\Phi (G,T)$
of roots of $G$ with respect to $T$. Then $(V,\Phi)$ forms a root
system.  For any $\alpha \in \Phi$ let $\alpha ^{\vee}$ denote the
corresponding coroot in $V^*$.  Let $\Delta=\{\alpha _1,
\dots ,\alpha _n \}$ be the set of simple roots corresponding to the choice
of $B$ and let $\lambda _1, \dots ,\lambda _n $ be the dual basis of
$\alpha _1^{\vee}, \dots , \alpha _n^{\vee}$.  This basis forms the set
of fundamental dominant weights with respect to $B$ and $\lambda _1,
\dots ,\lambda _n $ are vertices of the Weyl chamber $\frak c$
corresponding to the choice of $B$.

Let $E$ be a principal $G$-bundle on $X$ and let $E_P$ be a reduction
of structure group of $E$ to some parabolic subgroup $P$ of $G$
containing $B$. The parabolic $P$ corresponds to some subset
$\{\alpha _i\} _{i\in I}$ of $\Delta$ (see 1.1), where $I\subset \{1,
\dots ,n\}$.  The facet of $\frak c$ with vertices $\{\lambda _i\}
_{i\in I}$ corresponds to $P$ and it will be denoted by the same
letter.

Let $K$ denote the function field of $X$.  Take a maximal torus
$T_K\subset E(G) _K$, which is contained in $E_P(P) _K$. If $T_K$ is not
split then let us pass to a separable cover $\pi \: Y\to X$, where
$T_{K'}=\pi ^*T_K$ splits ($K'$ is the function field of $Y$). Let
$B_{K'}$ be the Borel subgroup contained in $\pi ^*E_P(P)_{K'}$.  We can
choose an isomorphism of root systems $\Psi = \Phi(\pi^*E(G)_{K'},T_{K'})$
and $\Phi $ such that the Weyl chambers corresponding to $B_{K'}$
and $B$ become equal and the facet corresponding to the parabolic $\pi
^* E_P(P) _{K'}$  corresponds to $P$. Let us set 
$$d(\frak c)=\sum _{i=1}^n\deg B_{K'} (V_{S(\alpha _i)})\cdot
\lambda_i ^{\vee},$$
where $V_{S(\alpha _i)}$ are considered with respect to $B$.  The
vector $d(\frak c)$ is well defined since $B_{K'}$ uniquely extends to
the Borel subgroup scheme of $\pi ^*E(G)$. It corresponds to the
vector $d(\frak c_{K'})$ that belongs to the complementary polyhedron
for $\Psi$ (see [Be2], Proposition 6.6).

\medskip 
{\sl 4.2.} Now let us consider $L$-modules $V_S$ defined with respect to the Levi
subgroup $L$ of $P$ (see 1.1).  All $V_S$ corresponding to shapes $S$
of level $1$ are called {\sl elementary $L$-modules}. The
corresponding vector bundles $E_L(V_S)$ are called {\sl elementary
vector bundles} (see [Be2], Definition 5.5) associated to $E_P$. The
degree of $E_L(V_S)$, where $S$ is a shape of level $1$ corresponding
to a root $\alpha$, is called the {\sl numerical invariant} of $E_P$
with respect to $\alpha $ and it is denoted by $n(E_P; \alpha)$. To be
compatible with [Be2] we need to define numerical invariants using
fundamental weights. Namely, let $\lambda$ be the fundamental weight
of $G$ with respect to $B$ dual to the coroot $\alpha ^{\vee}$. Then
we set $n(E_P; \lambda )=n(E_P; \alpha)$, defining numerical
invariants with respect to fundamental weights corresponding to $P$.

One can see that the degree of $E_L(V_S)$ for a shape $S$ of positive level
is a non-negative linear combination of numerical invariants of $E_P$
with respect to simple roots $\alpha \in I$, where $I$ is the subset
of the set of simple roots corresponding to $P$ (see Lemma 4.5).

\proclaim{Lemma 4.3}  
Let $\pi \: Y\to X$ be a finite morphism of normal projective 
varieties over $k$.  Let $E_P$ be the canonical reduction of a 
rational $G$-bundle $E$ defined over $X$. Then the numerical 
invariants of $\pi ^*E_P$ are non-negative. 
\endproclaim 
 
\pr In case $\pi$ is separable the lemma follows from [Be2], Lemma 7.1. 
So we can assume that $\pi$ is equal to the Frobenius morphism $F\:X\to X.$ 
In this case the lemma follows immediately from the fact that 
$(F^*E_P)(V_S)=F^*(E_P(V_S))$ and the degree of a rational vector bundle 
multiplies by $p=\char k$ under the Frobenius morphism, Q.E.D. 
 
\medskip

{\sl 4.4.} Now let us set
$$U(P)=\{\alpha \in \Phi\: \hbox{ there exists }\lambda \in \vert P
\hbox{ such that }\langle \alpha , \lambda ^{\vee}\rangle >0\}.$$
This set is equal to $\Phi ^+-\Phi _I$ and it has a natural
decomposition into subsets corresponding to roots of the same level
(with respect to $P$).  In the above notation $U(P)$ decomposes into
subsets
$$\Psi (P, \sum _{i\in I} n_i\lambda_i)=\{\alpha \in \Phi\: \langle
\alpha , \lambda _i ^{\vee}\rangle =n_i\hbox{ for every }i\in I \} .$$
Note that $\Psi (P, \sum _{i\in I} n_i\lambda_i)$ is the set of all
roots of shape $S=\sum _{i\in I} n_i\alpha _i$. Moreover, we have 
$$\rk  E_P (V_S)=|\Psi (P, \sum _{i\in I} n_i\lambda_i)|$$
and
$$\deg E_P (V_S)=\sum _{\alpha \in \Psi (P, \sum _{i\in I}
n_i\lambda_i)}\langle \alpha, d(\frak c)\rangle .$$

\proclaim{Lemma 4.5}
Let $\alpha=\sum _{i=1}^n n_i\alpha _i\in U(P)$. Then
$$\mu (E_P(V_{S(\alpha)}))= \sum _{i\in I} n_i {n(P, \lambda _i)\over
|\Psi (P, \lambda _i)|}=\sum _{i\in I} n_i \mu (E_P (V_{S(\alpha
_i)})).$$
\endproclaim

\pr
Let us set 
$$y(P)=\sum _{\lambda \in \vert P}{ n(P; \lambda)\over |\Psi (P, 
\lambda)|}\cdot {\lambda ^{\vee}}.$$
Then for all vertices of $P$ we have
$\langle \lambda, d(\frak c)\rangle =\langle \lambda, y(P)\rangle .$
Now let us note that $\sum _{\alpha \in \Psi (P, \sum _{i\in I}
n_i\lambda_i)} \alpha$ belongs to the vector space spanned by
$\{\lambda _i\} _{i\in I}$ (cf.~[Be2], Lemma 3.6). Hence
$$\deg E_P (V_{S(\alpha )})= \sum _{\alpha \in \Psi (P, \sum _{i\in I}
n_i\lambda_i)}\langle \alpha, y(P)\rangle =|\Psi (P, \sum _{i\in I}
n_i\lambda_i)|\cdot \sum _{i\in I}n_i {n(P, \lambda _i)\over |\Psi (P,
\lambda _i)|} ,$$
Q.E.D.

\heading 5. Asymptotic strong Harder--Narasimhan filtration\endheading

In the vector bundle case we proved that some Frobenius pull back has
strong Harder--Narasimhan filtration (see [La1], Theorem 2.7). We
prove that the same result holds for principal $G$-bundles.

\medskip 

\proclaim{Theorem 5.1} 
Let $E$ be a principal $G$-bundle on $X$ defined over a field $k$ of 
positive characteristic $p$ and let $F\: X\to X$ be the Frobenius 
morphism.  Then there exists $l$ such that $(F^l)^*E$ has strong 
Harder--Narasimhan filtration. 
\endproclaim

\pr 
To each rational $G$-bundle $E_P$ one can associate a sequence of 
parabolic subgroups $B\subset P_k$ corresponding to canonical 
reductions $E_{k,P_k}$ of $(F^k)^*E$. Since elements of this sequence 
are chosen from a finite set of parabolic subgroups that contain $B$, 
we can take a constant subsequence $\{ P_{j_k}\}$ of $\{P_k\}$. Set 
$P=P_{j_k}$ and let $I=\{\alpha _1, \dots , \alpha _r\}\subset \Delta$ 
be the corresponding subset determining $P$.  Then we consider $r$ 
sequences $\{n_{1, k}\}, \dots , \{n_{r, k}\}$ defined by 
$$n_{i, k}= {n (E_{j_k, P_{j_k}}; \alpha _i) \over p^{j_k}}$$ for
$i=1, \dots, r$ (see 4.2). By Lemma 4.5 there exist some nonnegative
rational numbers $A_1, \dots , A_r$, depending only on the type of $G$
and $P$, such that
$$\sum A_i n_{i, k}={\deg E_{j_k, P_{j_k}}\over p^{j_k}}.$$ Since the 
sequence $\left\{ {\deg E_{j_k, P_{j_k}}\over p^{j_k}} \right\}$ converges to a 
finite number $\deg _{HN, \infty}E$ we can find a subsequence $\{ 
j_{k_l}\}$ of $\{ j_k\}$ such that all the sequences $\{ n _{i, j_{k_l}}\}$ 
converge.  
 
Then for any  $\epsilon >0$ we can find $j_{k_{l}}$ and $j_{k_{l'}}$, 
where $l'>l$, such that 
$${n_{i,j_{k_{l'}}}\le {(1+\epsilon)} \cdot n _{i, j_{k_l}}}$$ 
for  $i=1, \dots ,r$. 
Let us replace $E$ with $(F^{j_{k_l}})^*E$ and set 
$k=j_{k_{l'}}-j_{k_{l}}$.  Let $E_P$ be the Harder--Narasimhan 
filtration of $E$ and assume that $(F^k)^*E_P$ is not the 
Harder--Narasimhan filtration of $(F^k)^*E$. Let ${\tilde E}_{P}$ be 
the Harder--Narasimhan filtration of $(F^k)^*E$. 
 
Let us set $\G =((F^k)^*E )(G)$, $\P= {\tilde E}_P (P)$, $\Q
=(F^k)^*(E_P(P))$ and let $K$ denote the function field of $X$. Then
there exists a maximal torus $T_K\subset \G _K$ which is contained in
$\P _K\cap \Q_K$ (see [SGA3], exp. XXVI, 4.1.1). If $T_K$ is not split
then let us pass to a separable cover $Y\to X$, where $T_K$ splits. By
Lemma 2.9 we may actually replace $X$ with $Y$ and assume that $T_K$
splits on $X$. Our assumption imply that $\P_K\ne \Q _K$. Let $B_K$
and $B_K'$ be the Borel subgroups contained in $\P _K$ and $\Q _K$,
respectively.
 
We can find such an isomorphism of root systems $\Phi (G, T)$ and
$\Phi (\G _K, T_K)$ that the Weyl chambers $\frak c$ and $\frak c _K$
corresponding to $B$ and $B_K$ become equal and the facet $P_K$
corresponding to $\P _K$ corresponds to the facet $P$ of $\frak c$. We
can also find such an isomorphism of root systems $\Phi (G, T)$ and
$\Phi (\G _K, T_K)$ that the Weyl chambers $\frak c$ and $\frak c _K$
corresponding to $B$ and $B_K'$ become equal and the facet $Q_K$
corresponding to $\Q _K$ corresponds to the facet $P$ of $\frak c$.
This shows that there exists an isomorphism $\sigma$ of $V=X^*(T_K)\ot
_{\Z} \R $ preserving $\Phi=\Phi (\G _K, T_K)$ and such that the image
of the facet $P_K$ is equal to the facet $Q_K$. Obviously, $\sigma$ is
an element of the Weyl group $W$.
 
Let $(\cdot \, , \cdot)$ be a scalar product on $V$ that gives a
$W$-invariant Euclidean metric on $V$. Let $\| \cdot \|$ be the
associated norm. This can be used to identify $V$ and $V^*$.  In this
identification if $\alpha \in \Phi$ then $\alpha ^{\vee}$ is
identified with $2\alpha \over (\alpha, \alpha)$.
 
Let $\frak C$ be the set of Weyl chambers of $\Phi$ and let
$d=(d(\frak c))_{\frak c\in \frak C} $ consists of vectors in $V^*$
that correspond to ${1\over p^k}$ of vectors in the complementary
polyhedron for $\Phi $ determined by $T_K\subset \G _K$ (see [Be2],
Proposition 6.6). This set still forms a complementary polyhedron such
that
$$n(P_K, \lambda )= {n({\tilde E}_P; \lambda)\over p^k}\ge 0$$ 
and 
$$n(Q_K, \lambda )= {n((F^k)^*E_P; \lambda)\over p^k}= n(E_P; \lambda )\ge 0$$ 
(cf. Lemma 4.3).  
 
Our assumptions show that $\sigma$ is an isometry such that $\sigma (P_K)=Q_K$ and  
$$n( P_K; {\lambda })\le (1+\epsilon)\cdot n (Q_K;\sigma (\lambda )) 
\leqno (*)$$ for all vertices $\lambda$ of $P_K$.  Let $\Psi (P, 
\lambda)$ be the elementary set of roots associated to facet $P$ and 
its vertex $\lambda$ (see 4.4). Obviously, $|\Psi (P_K;\lambda
)|=|\Psi (Q_K;\sigma (\lambda ))|$ for any vertex $\lambda $ of $P_K$.
 
As in the proof of Lemma 4.5   let us set  
$$y(P)=\sum _{\lambda \in \vert P}{n(P; \lambda)\over |\Psi (P, 
\lambda)|}\cdot {\lambda ^{\vee}}$$
for a facet $P$.  By assumption $y(P_K)\in P_K^{\vee}$ and $y(Q_K)\in
Q_K^{\vee}$ (this follows from non-negativity of the numerical
invariants of $P_K$ and $Q_K$).
 
If we consider $y(Q_K)$ and $y(P_K)$ as vectors in $V$, then [Be2], Lemma 2.5 and
[Be2], Proposition 3.13 show that 
$$( \lambda, y(Q_K))=\langle \lambda, d(\frak c_K')\rangle \le ( \lambda, y(P_K))$$  
for any vertex $\lambda$ of $Q_K$.
This implies that 
$$(y(Q_K), y(Q_K))\le (y(Q_K), y(P_K)).\leqno (**)$$ Since $(y(Q_K), 
y(P_K))=\| y(Q_K)\| \cdot \| y(P_K)\| \cdot \cos \alpha$, where $\alpha $ is 
the angle between $y(Q_K)\in Q_K^{\vee}$ and $y(P_K)^{\vee}$. Since 
$\P \ne \Q$, this angle is non-zero and $\cos \alpha$ cannot be larger 
than the maximum $s_0$ of $\cos \alpha$ over all non-zero angles 
between different facets in a partition of $(X^*(T)\ot_{\Z} \R )^*$ 
determined by $\Phi (G ,T)^{\vee}$. Obviously, $s_0<1$, since $\Phi (G 
,T)$ is finite. Hence by $(**)$  
$$\| y(Q_K)\| \le s_0 \| y(P_K)\| .$$ 
But using $(*)$ we get  
$$\aligned 
\| y(P_K)\| ^2 &=\sum _{\lambda , \mu \in \vert P_K} {  n(P_K; \lambda)
\cdot n(P_K; \mu )\over |\Psi (P_K, \lambda)|  
\cdot |\Psi (P_K, \mu)|}\cdot  (\lambda ^{\vee}, \mu ^{\vee})\\ 
&\le (1+\epsilon )^2 \sum _{\lambda , \mu \in \vert P_K} {
n(Q_K;\sigma(\lambda))\cdot n(Q_K; \sigma(\mu ))
\over |\Psi (Q_K, \sigma(\lambda ))| \cdot |\Psi (Q_K, \sigma(\mu ))|}
\cdot (\sigma(\lambda ^{\vee} ),\sigma (\mu ^{\vee}))\\ 
&=(1+\epsilon)^2 \|x (Q_K) \| ^2.\\
\endaligned 
$$ 
Therefore $1\le s_0 (1+\epsilon)$ and for small $\epsilon$ we get a
contradiction, Q.E.D.
 
\bigskip

\proclaim{Corollary 5.2} 
Let $E$ be a principal $G$-bundle on a curve $C$ defined over a field
$k$ of positive characteristic $p$ and let $F\: C\to C$ be the
Frobenius morphism.  Then there exists $l$ such that the
Harder--Narasimhan filtration ${\tilde E}_P$ of $(F^l)^*E$ has a
strongly semistable reduction to the Levi component $L\subset P$.
\endproclaim 

\pr
Let us take $l$ such that $(F^l)^*E$ has the strong Harder--Narasimhan
filtration. We can take even larger $l$ such that the degrees of all
elementary vector bundles are greater than $\deg K_C$. The latter can be
easily achieved, since the degrees of elementary vector bundles are
positive and they multiply by $p$ under the Frobenius pull back.

The theorem will be proved if we show that the non-abelian
cohomology group $H^1(C, E_P (R_u(P)))$ is trivial (cf. [SGA3],
exp. XXVI, Corollaire 2.2 and [Su], Theorem 2.2).

Now note that $E_P(R_u(P))$ has a filtration whose quotients are
vector bundles $E_P(V_S)$. Since the non-abelian cohomology of a
vector bundle coincides with its usual sheaf cohomology, it is
sufficient to prove that the sheaf cohomology $H^1(C, E_P(V_S))$
vanish. But by Serre duality $h^1(C, E_P (V_S))=h^0 (C, E_P(V_S)^*\ot
\omega _C)$ and $ E_P(V_S)^*\ot \omega _C$ is a semistable vector bundle 
of negative degree (semistability follows from Theorem 2.7), so it has
no sections, Q.E.D.

\medskip

The above corollary generalize [Ra1], Lemma 3.7 (where $C=\P ^1$)
and [Su], Theorem 2.2 (where $C$ is an elliptic curve). In the vector
bundle case it was first proved by V.~B.~Mehta and S.~Subramanian and
the author learned it from S.~Subramanian. In this case the corollary
says that if $E$ is a vector bundle on a curve then the
Harder--Narasimhan filtration of some Frobenius pull back of $E$
splits into a direct sum.

\heading 6. Semistability of Frobenius pull backs \endheading

\medskip
{\sl 6.1.} Let us fix a maximal torus $T$ and a Borel subgroup
$B\supset T$ in $G$. Let $E$ be a principal $G$-bundle with canonical
reduction $E_P$ for some parabolic subgroup $P$ containing $B$. 
Let $\lambda$ be the fundamental weight corresponding to one of
vertices of the facet corresponding to $P$. Let $V_\lambda $
denote the elementary module $V_{S(\alpha)}$, where $\alpha$ is the
simple root such that $\lambda$ is dual to the coroot $\alpha ^{\vee}$. 

By abuse of notation we will use the same notation to denote an
elementary module corresponding to other parabolic subgroups of
$G$. This will not lead into confusion since we use this notation
together with a reduction of structure group of $E$ to the parabolic
with respect to which we consider it as an elementary module.

Let us recall that a parabolic subgroup of $G$ is called {\sl
maximal}, if it is proper (i.e., different to $G$) and if it is not
contained in any other proper parabolic subgroup of $G$.
In the vector bundle case it is obvious that any component of the
Harder--Narasimhan filtration destabilizes the vector bundle. The same
fact holds for principal $G$-bundles:

\proclaim{Proposition 6.2}
Let $Q\subset G$ be a maximal parabolic subgroup containing $P$ and
let $E_Q$ be the extension of structure group of $E_P$ to $Q$. 
Let $\mu $ be the fundamental weight corresponding to $Q$.
Then
$$\mu ( E_Q(V_{\mu}))= \sum _{\lambda \in \vert P} {\langle \mu,
\lambda ^{\vee}\rangle \over \langle \mu, \mu^{\vee}\rangle } \mu 
(E_P(V_{\lambda})),\leqno (6.2.1)$$ 
where $\vert P$ denote the set of fundamental weights corresponding to $P$.
In particular, the degree of $E_Q (Q)$ is positive.
\endproclaim

\pr  
Let $\frak c$ be the Weyl chamber corresponding to $B$ and let
$d(\frak c)$ be as in 4.1. This vector will be used to compute the
degree of $E_Q(Q)$.  Let $P'$ and $Q'$ be the facets corresponding to
$E_P(P)$ and $E_Q(Q)$, respectively.  The facet $Q'$ has only one
vertex $\mu$ and it is also one of the vertices of $P'$.

Now let us note that
$$\langle \mu , d(\frak c)\rangle =\langle \mu , y(P')\rangle= \sum
_{\lambda\in \vert P'}{n(P', \lambda)\over |\Psi (P',\lambda)|}\langle
\mu , {\lambda}^{\vee} \rangle = \sum _{\lambda \in \vert P} \langle \mu , 
{\lambda}^{\vee} \rangle \cdot {\mu (E_P(V_{\lambda}))}   ,$$
since $\mu$ is a vertex of $P'$. But we also have
$$\langle \mu , d(\frak c)\rangle =\langle \mu , y(Q')\rangle= 
{n(Q',\mu )\over |\Psi (Q', \mu)|}\langle \mu , {\mu}^{\vee}\rangle =
\langle \mu , {\mu}^{\vee}\rangle \cdot {\mu ( E_Q(V_{\mu}))} ,$$
since $\mu$ is also a vertex of $Q'$. Comparing these two equalities
yields the required equality.

Now note that the degree of $E_Q(Q)$ is a positive multiple of $\mu (
E_Q(V_{\mu}))$ (this follows, e.g., from Lemma 4.5).  Since
$E_P$ is the canonical reduction we have $\mu ( E_P(V_{\lambda}))>0$,
so the inequality $\deg E_Q(Q)>0$ follows from the fact that the
coefficients in (6.2.1) are non-negative and one of them is equal to
$1$, Q.E.D.

\medskip

The next theorem bounds the slopes of elementary vector bundles of
Frobenius pull back of a principal $G$-bundle. The method of proof is
similar to the proof of [La1], Corollary 6.2.

\proclaim{Theorem 6.3} 
Let $E$ be a semistable rational $G$-bundle which is not strongly
semistable. Let us take $l$ such that $\tilde E=(F^l)^*E$ has 
strong Harder--Narasimhan filtration and let ${\tilde E}_P$ be the
canonical reduction of $\tilde E$. Let $\mu$ be a fundamental weight
of $P$ and let $Q$ be the corresponding maximal parabolic containing
$P$.  Then for some integer $0\le i<l$ we have
$$0<\mu ({\tilde E}_Q(V_\mu))\le \mu _{\max} ((F^{i})^*\Omega _X),$$
where ${\tilde E}_Q$ is the extension of structure group of ${\tilde
E}_P$.
\endproclaim

\pr 
Note that ${\tilde E} _Q$ does not descend $l$ times under the
Frobenius morphism, since by Proposition 6.2 it would contradict
semistability of $E$. Let $i$ denote the non-negative integer such
that ${\tilde E} _Q$ descends $i$ times under Frobenius, but it does
not descend $(i+1)$ times. Let us write ${\tilde E}_Q=(F^i)^*{E}_Q$
for some reduction $E_Q$ of $(F^{l-i})^*E$.

Let $\sigma \: X\to (F^{l-i})^*E/Q $ be the section corresponding to
the reduction $E_Q$.  Then we have a map $\tau \: T_X\to
\sigma ^* N_{\sigma}$, where $N_\sigma$ is the normal bundle of $\sigma
(X)$ in $(F^{l-i})^*E/Q$ (see, e.g., [MS], the proof of Theorem 4.1).
Moreover, $\sigma^*N_{\sigma}= E_Q(\frak g/\frak q)$ and the map
$\tau$ is non-zero, because otherwise $\sigma$ would descend under the
Frobenius morphism, contradicting our assumption on $i$.  Hence the
map $(F^i)^*(T_X)\to {\tilde E}_Q(\frak g/\frak q)$ is also
non-zero. This implies that
$$\mu _{\min} ((F^i)^*(T_X))\le \mu _{\max}({\tilde E}_Q(\frak g/\frak q)).$$ 
But $\mu _{\max}({\tilde E}_Q(\frak g/\frak q))=-\mu _{\min}( \tilde E_Q( \frak u)),$ where
$\frak u$ is the Lie algebra of the unipotent radical $R_u(Q)$ of $Q$.
Hence we get
$$\mu _{\min}(\tilde E_Q(\frak u))\le \mu _\max ((F^i)^*\Omega _X).$$
Note that $\tilde E_Q(\frak u)$ has a filtration with strongly
semistable quotients, whose slopes are equal to multiples of $\mu
(\tilde E_Q(V_{\mu}))$ (by Lemma 4.5) and $\tilde E_Q(V_\mu)$ is a
quotient of $\tilde E_Q(\frak u)$. Therefore $\mu _{\min}(\tilde
E_Q(\frak u))=\mu (\tilde E_Q(V_\mu ))$, Q.E.D.

\proclaim{Corollary 6.4} {\rm ([MS], Theorem 4.1)}
If $\mu _{\max}(\Omega _X)\le 0$ then every semistable rational
$G$-bundle is strongly semistable. In particular, every rational
$G$-bundle on $X$ has strong Harder--Narasimhan filtration.
\endproclaim

\pr 
If there exists a semistable $G$-bundle which is not strongly
semistable then by Theorem 5.1 there also exists a semistable
$G$-bundle $E$ such that $F^*E$ is not semistable, but it has 
strong Harder--Narasimhan filtration. But this contradicts Theorem
6.3, Q.E.D.

\medskip
{\sl 6.5.} Let $P$ be a parabolic subgroup of $G$. Then by [Be2], Proposition 
1.9 there exist some positive numbers $b_{\mu, P}$ such that
$$\sum _{\alpha \in  U(P)}\alpha =\sum _{\mu\in \vert P}b_{\mu, P}\mu.$$

\proclaim{Corollary 6.6}
Assume that $\mu _{\max}(\Omega _X)> 0$. Let $E$ be a semistable
$G$-bundle and let $P$ be the parabolic subgroup of $G$ corresponding
to the strong Harder--Narasimhan filtration of some Frobenius pull
back of $E$. Then
$$\deg _{\HN ,\infty} {E}\le \left(\sum _{\mu \in \vert P}b_{\mu, P}\langle \mu,
\mu ^{\vee}\rangle \right) {L _{\max} (\Omega _X)\over p}.$$
In particular, we have
$$\deg _{\HN}E(\frak g)\le \deg _{\HN ,\infty } E(\frak g)\le
{2\dim\frak g\over p}\cdot \left(\sum _{\mu \in \vert P}b_{\mu,
P}\langle \mu,
\mu ^{\vee}\rangle \right) {L _{\max} (\Omega _X)}.$$
\endproclaim

\pr Let us take $l$ such that both $(F^l)^*E$ and $(F^{l-1})^*(\Omega _X)$
have strong Harder--Narasimhan filtrations. Then $\mu _{\max}
((F^{l-1})^*\Omega _X)=p^{l-1} L _{\max} (\Omega _X)$.
Note that
$$\aligned
\deg _{\HN ,\infty} {E}&=\langle \sum _{\alpha \in  U(P)}\alpha , y(P)\rangle =
\sum _{\lambda \in \vert P}\mu (E_P (V_{\lambda}))\, \langle \sum _{\alpha \in  U(P)}\alpha , 
\lambda^{\vee}\rangle \\
&=\sum _{\mu \in \vert P}b_{\mu ,P }\sum _{\lambda \in \vert P}\mu
(E_P (V_\lambda ))\, \langle \mu ,\lambda^{\vee}\rangle .\\
\endaligned
$$
Hence the first inequality follows from Proposition 6.2 and Theorem 6.3. The second 
inequality follows from Corollary 3.4, Q.E.D.

\medskip
Let us set 
$$b(G)=2\dim \frak g\cdot  \mathop{\max} _{P\subset G} 
\left(\sum _{\mu \in \vert P}b_{\mu, P}\langle \mu,
\mu ^{\vee}\rangle \right) ,$$
where the maximum is taken over a finite set of all parabolic subgroups of
$P$ containing $B$.

\proclaim{Corollary 6.7}
Assume that $\mu _{\max} (\Omega _X)> 0$ and $p>b(G)\cdot L_{\max}(\Omega
_X)$. Then $E$ is semistable if and only if $E(\frak g)$ is
semistable.
\endproclaim

\pr 
One implication follows from Corollary 2.8. The other implication
follows from Corollary 6.6 and the remark that $\deg _{\HN}E(\frak
g)$ is an integer, Q.E.D.

\medskip 
Corollary 6.7 is similar to, but usually weaker than, the main result
of [IMP] applied to the adjoint representation. However, Corollary 6.6
bounds the degree of instability of the adjoint bundle $E({\frak g})$
even in small characteristic, so in the cases when [IMP] gives no
information. In the same way one can use Theorem 8.4 to prove that an
extension of the structure group of a semistable $G$-bundle is
semistable if the characteristic $p$ is large. In the case when the
corresponding group homomorphism is a representation we get a weak
form of the main result of [IMP].

\heading 7. Boundedness of principal $G$-bundles\endheading

\medskip
{\sl 7.1.} An {\sl algebraic family} of rational $G$-bundles on $X$
parametrised by $S$ is a rational $G$-bundle on $X\times S$, whose
restriction to each fibre of the projection $X\times S\to S$ is a
rational $G$-bundle. A family $\E$ of principal $G$-bundles on $X$ is
called {\sl bounded} if there exists an algebraic family of principal
$G$-bundles on $X$ parametrised by a scheme of finite type over $k$
and such that it contains each element of $\E$ (up to an isomorphism).

\medskip

{\sl 7.2.}  Let us recall that the {\sl degree} of a principal
$G$-bundle $E$ is a homomorphism $d_E\: X^*(G)\to \Z$ given by $\chi
\to \deg E(\chi)$, where $E(\chi )$ is the line bundle associated to
$E$ by $\chi$.

Note that the character group $X^*(G)$ of $G$ is a subgroup of finite
index in $X^*(R(G))$ (this follows from the well known facts saying
that $G=R(G)\cdot (G,G)$ and $R(G)\cap (G,G)$ is finite). Hence we can
uniquely extend $d_E$ to a homomorphism $X^* (R(G))\to
{\Bbb Q}$.  This homomorphism should be thought of as a slope of $E$ (look
at the $G=\GL (V)$ case).

Since $R(G)\subset T$ we also get the induced homomorphism $X^*(T)\to
X^*(R(G))\to {\Bbb Q}$ denoted by $d'_E$. Note that 
$$X^*(G)=X^*(T)_0= \{ \lambda \in X^*(T) \: \, \langle \lambda ,
\alpha ^{\vee}\rangle =0 \hbox{ for all }\alpha \in \Phi\} $$
and the restrictiction of $d'_E\: X^*(T)\to {\Bbb Q}$ to $X^*(G)\subset X^*(T)$
induces the original degree homomorphism $d_E\: X^*(G)\to \Z$. 

We can interpret $d'_E$ in the following way. Let 
$$X^*(T)_+= \{ \lambda \in X^*(T) \: \, \langle \lambda ,
\alpha ^{\vee}\rangle \ge 0 \hbox{ for all }\alpha \in \Phi ^+\} $$
denote the dominant weights of $T$ with respect to the set of positive
roots $\Phi ^+$.  For any $\lambda \in X^*(T)_+$ we denote by
$L(\lambda)$ the simple $G$-module with highest weight $\lambda$.  The
radical $R(G)$ is contained in the centre $Z(G)$ and hence it acts on
$L(\lambda)$ through the restriction of $\lambda$ to $R(G)$ (see [Ja],
Part II, 2.10).  This implies that the radical $R(G)$ acts on $\det
E(L(\lambda ))$ through the restriction of $\lambda ^{\dim
L(\lambda)}$ to $R(G)$.  Hence
$$d'_E(\lambda)= \mu (E(L(\lambda))).$$

\medskip

\proclaim{Theorem 7.3}  
Let us fix some constant $C$. Then the family $\F$ of all semistable
principal $G$-bundles $E$ on $X$ such that the degree of $E$ is fixed
and $a_2(\t E(\frak g )) \ge C$ (see 1.4), is bounded.
\endproclaim

\pr  
By Corollary 6.6 and Proposition 2.4 the slopes of maximal
destabilizing subsheaves of $\t E(\frak g)$ for $E\in \F$ are bounded
from the above by some (explicit) constant $C'$. Then Theorem 1.5
implies that the family $\F '= \{E (\frak g)\} _{E\in\F}$ of rational
vector bundles is bounded.

Let $G'$ be the image of the adjoint representation $\Ad\: G\to \GL
(\frak g)$. It is equal to the quotient of $G$ by the centre group
scheme $Z(G)$.  Let $\E =\{E_{G'}\} _{E\in \F}$ be the family of
$G'$-bundles obtained from $G$-bundles in $\F$ by extension of
structure group to $G'$.  Since $G' \hookrightarrow \GL (\frak g)$,
each rational $G'$-bundle in $\E$ can be constructed as a reduction of
structure group of a rational $\GL (\frak g)$-principal bundle $F$
from $\F'$.  But such reductions of structure group can be
parametrized by a scheme of finite type corresponding to sections
$U\to F(\GL (\frak g)/G')$ defined on some big open subset $U$ of $X$.
Since the family $\F '$ is bounded, this shows that the family $\E$ is
also bounded (see [Ra3], Lemma 4.8.1 for a precise argument).

For any character $\lambda\in X^*(T)$ let $k_{\lambda}$ be a
one-dimensional $B$-module whose restriction to $T$ is $\lambda$.  Let
$H^0(\lambda)$ denote the $G$-module $H^0(G/B, \L _\lambda)$, where $\L
_{\lambda}$ is the line bundle associated to the $B$-module
$k_{\lambda}$. Note that $H^0(\lambda)$ are usually not simple
$G$-modules, but the centre group scheme $Z(G)$ acts on $H^0(\lambda)$
through scalars (see [Ja], Chapter II, Proposition 2.8). In fact, it
acts through the restriction of $\lambda$ to $Z(G)$.

Let $\rho \: G\to \GL (V)$ be any finite-dimensional faithful rational
representation, which is a direct sum of $G$-modules of the form
$H^0(\lambda)$ for some $\lambda\in X^*(T)$. Then for any direct summand
$H^0(\lambda)$ of $V$ we have a group homomorphism $G'\to \PGL
(H^0(\lambda))$ induced from $\rho$.  Hence the family $\F _\lambda$ of
$\PGL (H^0(\lambda))$-bundles obtained from the family $\F$ by extension
of structure group $G\to G'\to \PGL (H^0(\lambda))$ is also bounded.

Let $\E _\lambda$ be the family of vector bundles associated to $\GL
(H^0(\lambda))$-bundles obtained from $\F$ by extension of structure
group $G\to \GL (H^0(\lambda))$. Two $H^0(\lambda)$-vector bundles give
the same principal $\PGL (H^0(\lambda))$-bundle if and only if their
projectivisations are isomorphic, i.e., they differ by tensoring by a
line bundle. Since the degree in the family $\F$ is fixed the degree
of vector bundles in the family $\E _\lambda$ is also fixed.  This and
the fact that the family $\F_\lambda$ is bounded imply that the family
$\E _\lambda$ is also bounded (see [Ra3], 4.15 for a general argument).

Hence the family $\E =\{E(V)\} _{E\in \F}$ of vector bundles, obtained
as direct sums of vector bundles from families $\E _\lambda$, is also
bounded (see [Ra3], the proof of Proposition 4.12). But, as before,
this implies that the family $\F$ is also bounded, Q.E.D.

\medskip
If in the above proof instead of Corollary 6.6 we use Corollary 8.5
(and Corollaries 2.8 and 6.4) then we get the following theorem:

\proclaim{Theorem 7.4} 
Let $C_1$ and $C_2$ be fixed constants. Then the family of all
rational $G$-bundles $E$ on $X$ such that the degree of the canonical
parabolic of $E$ is $\le C_1$, the degree of $E$ is fixed, and $a_2(\t
E(\frak g))\ge C_2$ is bounded.
\endproclaim

As a special case of the above theorem we get Theorem 0.1.

\heading 8. Instability of bundles associated to representations\endheading

In this section we give an explicit bound for $L_{\max}(E({\frak h}))$
for a semistable $G$-bundle $E$ and a homomorphism $\rho :G\to H $ of
semisimple groups (see Theorem 8.4). In particular, we get the bound
for the difference $\mu _{\max}(E(V))-\mu _{\min}(E(V))\le
L_{\max}(E(\gl V))$ for an arbitrary representation $G\to \GL (V)$ of
a semisimple group.

\medskip

{\sl 8.1.} Let us first start with the simplest case of $G=\SL (V)$,
where the bound is particularly strong.

Set $n=\dim V$ and let $T_{\GL (V)}$ be the standard maximal torus of
$\GL(V)$. Let $\epsilon _i$, $ i=1,\dots , n$ be the standard basis of
$X^*(T_{\GL (V)})$, i.e., $\epsilon _i$ is the restriction of the
matrix coefficient $x_{ii}$ to $T_{\GL (V)}$. Let $T_{\SL
(V)}=T_{\GL(V)}\cap
\SL(V)$ be the corresponding maximal torus in $\SL (V)$.  Then
$X^*(T_{\SL (V)})=X^*(T_{\GL (V)})/\Z (\epsilon _1+\dots +\epsilon
_n)$. Let $\omega_i=\epsilon_1+\dots+\epsilon _i$ for $ i=1,\dots , n$
be the dominant weights of $\GL (V)$. Then the restrictions
$\omega_i'=\omega _i|_{T_{\SL(V)}}$ for $i=1,\dots , n-1$ are the
dominant weights of $\SL (V)$.

Let $W$ be a polynomial $\GL (V)$-module. Then there exists $m$ such
that $W$ is a submodule of $V^{\ot m}$ (see [KP], Proposition 5.3). In
particular, if $W=L(\lambda)$ is the simple $\GL (V)$-module with
highest weight $\lambda=\sum _{i=1}^n m_i \omega _i$, then $m$ is
uniquely determined by $\lambda$ and it is equal to the {\sl degree}
$|\lambda|=\sum _{i=1}^{n} im_i$. This follows from the facts that the
scalars act on $L(\lambda)$ through the restriction of $\lambda$ and
$L(\omega _i)=\wedge ^i V$.

Every dominant weight of $\SL (V)$ can be written as a sum
$\lambda=\sum_{i=1}^{n-1} m_i\omega_i'$ of fundamental weights.  Then
the corresponding weight $\lambda'=\sum_{i=1}^{n-1} m_i\omega_i$ of
$\GL (V)$ is polynomial (see, e.g., [Ja], Proposition A.3). Hence by
the above the $\GL (V)$-module $L(\lambda ')$ is a submodule of
$V^{|\lambda '|}$. But $L(\lambda ')$ is the simple $\SL (V)$-module
with highest weight $\lambda=\lambda'|_{T_{\SL (V )}}$ (see [Ja], II,
2.10.(2)).

Hence the simple $\SL (V)$-module $L(\lambda)$ with highest weight
$\lambda$ is an $\SL (V)$-submodule of $ V^{\otimes |\lambda|}$, where
$|\lambda|=\sum _{i=1}^{n-1}im_i$ is the {\sl degree} of $\lambda$.

Let $W$ be an $\SL (V)$-module. Then the maximum of degrees of
fundamental weights, whose modules occur as quotients of the
Jordan--H\"older filtration of $W$ is called the {\sl $\JH$-degree} of
$W$ and denoted by $\JH (W)$.

\proclaim{Lemma 8.2}
Let $W$ be an $\SL (V)$-module. Let $E$ be a principal $\SL
(V)$-bundle and let $E(V)$ and $E(W)$ denote the associated vector
bundles. Then
$$\JH (W)\cdot L_{\min} (E(V))\le L_{\min} (E(W))\le 
L_{\max} (E(W))\le \JH (W)\cdot L_{\max} (E(V)).$$ 
\endproclaim

\pr 
If $W$ is not a simple $\SL (V)$-module then take the Jordan--H\"older
filtration of $V$ and let $V_i$ denote the quotients of this
filtration.  Since $L_{\max} (E(W))$ is less or equal to the maximum
of $L_{\max}(E(V_i))$ and $V_i$ are simple, we can assume that $W$ is
simple. Then $W$ is isomorphic to some $L(\lambda)$.  Since
$L(\lambda)$ is a submodule of $V^{\ot |\lambda|}$, so $E(W)$ is a
subbundle of $E(V)^{\otimes |\lambda|}$.  Therefore by 1.2
we have
$$L_{\max}(E(W))\le |\lambda|\cdot L_{\max}(E(V)).$$ The second
inequality follows from the above one applied to the dual
representation, Q.E.D.

\medskip

\proclaim{Corollary 8.3}
Let $\rho \: \SL (V)\to \GL (W)$ be a homomorphism and let $E$ be a
principal $\SL (V)$-bundle. Then
$$L_{\max} (E(\gl W))=L_{\max}(E(W))-L_{\min}(E(W))\le \JH (\rho
')\cdot L_{\max} (E(V)),$$ where $\rho' =\Ad _{\GL (W)}\circ \rho$.  

In particular, if $E$ is  semistable, $\char k=p$
and $\mu _{\max}(\Omega _X)>0$ then
$$L_{\max}(E(\gl W))\le (\dim V-1)\JH (\rho')\cdot {L_{\max}(\Omega
_X)\over p}.$$
\endproclaim

\pr
The first inequality follows from Lemma 8.2 applied to the
representation $\rho '$. The second inequality follows from the first
one and Theorem 1.3, Q.E.D.

\medskip
We can apply a similar method to prove a theorem similar to the
second part of Corollary 8.3 for any homomorphism of semisimple groups
(in fact, the statement is slightly more general):

\proclaim{Theorem 8.4}
Let $\rho \: G\to H$ be a homomorphism of connected reductive groups
over $k$.  Assume that $\rho(R(G))\subset R(H)$. Let $E_G$ be a
semistable principal $G$-bundle and let $E_H$ be the extension of
structure group of $E_G$ to $H$.
\item{(1)} If $\char k=0$ or $\mu _{\max} (\Omega _X)\le 0$ then $E_H$ is 
strongly semistable.
\item{(2)} If $\char k=p$ and $\mu _{\max} (\Omega _X)> 0$ then there 
exists some explicit constant $C(\rho)$ depending only on $\rho$ such that
$$0\le L _{\max}(E_H(\frak h))\le C(\rho )\cdot {L_{\max} (\Omega _X)\over
p} .$$ In particular, if $p$ is large then both $E_H(\frak h)$ and
$E_H$ are semistable.
\endproclaim

\medskip

\pr 
(1) follows from Theorem 2.7 and Corollary 6.4. Hence we can assume
that we are in case (2).

Let $\lambda\: G_m\to G$ be a $1$-parameter subgroup of $G$.  Then we
can associate to $\lambda$ a closed subgroup $P(\lambda)$ of $G$ by
$$P(\lambda )=\{p\in G\: \, \lim _{t\to 0}\lambda(t)\cdot p \cdot \lambda
(t)^{-1} \hbox{ exists in }G\} .$$ It is a parabolic subgroup and any
parabolic subgroup $P$ of a reductive group $G$ is of this form for
some $1$-parameter subgroup $\lambda$ (see [Sp], Proposition 8.4.5).
The unipotent radical $R_u(P(\lambda))$ of $P(\lambda)$ consists of
such points $p\in P(\lambda)$ that $\lim _{t\to
0}\lambda(t)\cdot p \cdot \lambda (t)^{-1}=e$, where $e$ is the
neutral element in $G$.

By Theorem 5.1 we can take such $l$ that ${\t E}=(F^l)^*E$ has the strong
canonical reduction ${\t E}_P$. Let $\lambda$ be a $1$-parameter
subgroup such that $P$ is associated to $\lambda$ and let $Q$ be the
parabolic subgroup of $H$ associated to $\rho \circ \lambda$. Then
$\rho (P)\subset Q$ and $\rho (R_u(P))\subset R_u(Q)$. There exists a
filtration of $\frak h$ with simple $Q$-modules as quotients and such
that $R_u (Q)$ acts trivially on each factor.  This can be constructed
by taking $\frak u\subset \frak q\subset \frak h$, where $\frak u$ is
the Lie algebra of $R_u(Q)$, and taking the corresponding filtrations
of $\frak u$, $\frak q/\frak u$ and $\frak h/\frak q$ (see 1.1). Now
take a further refinement of this filtration $V_m\subset V_{m-1}
\subset \dots \subset V_0=\frak h$ such that the quotients are simple
$P$-modules. By construction $R_u (P)$ acts trivially on each quotient
$W_i=V_i/V_{i+1}$ of this filtration and hence ${\t E}_H(V_i)/{\t
E}_H(V_{i+1})={\t E}_L(W_i)$. Since $W_i$ is a simple $L$-module, by
Schur's lemma the radical of $L$ (which is contained in $Z(L)$) acts
on $W_i$ by scalars. In particular, the above filtration gives rise to
a filtration of ${\t E}_H(\frak h)$ with strongly semistable
quotients. Degrees of these quotients can be determined in the
following way.  Note that $P\to P/R_u(P)=L$ induces the map $X^*(L)\to
X^*(P)$ of character groups, which we compose with the degree map
$d_{{\t E}_P}\: X^*(P)\to {\Z}$. In this way we get the degree map
$d_{{\t E}_L}:X^*(L)\to \Z$. As in 7.2 we can extend it to $d_{{\t
E}_L}:X^*(T_L)\to {\Bbb Q}$, where $T_L$ is a maximal torus in $L$.
If $W_i$ is a simple $L$-module with highest weight $\lambda _i\in
X^*(T_L)$, then the slope of ${\t E}_L(W_i)$ can be computed as $d_{\t
E_L}(\lambda_i)$. Writing $\lambda _i$ as a sum of fundamental weights
of $L$, we can use Proposition 6.2 and Theorem 6.3 (or Corollary 6.6),
to bound $d_{\t E_L}(\lambda_i)$ by means of the coefficients in the sum
times $p^{l-1}L _{\max}(\Omega _X)$. In particular, since 
$$L_\max (E(\frak h))={\mu_{\max}({\t E}({\frak h}))\over p^l}\le 
\mathop{\max}_i\, {\mu ({\t E}_L(W_i))\over p^l},$$ 
this gives the required explicit bound on $L_{\max}(E(\frak h))$, 
Q.E.D.

\medskip
{\sl Remarks.} 

(1) Note that the above theorem also bounds $\deg _{\HN}E_H$. This
follows from the definition, since $\deg _{\HN}E_H$ is the degree of a
subbundle of the degree zero vector bundle $E_H(\frak h)$.

(2) From the proof of the above theorem one can easily see that
$C(\rho)$ can be explicitly bounded by means of the heights of the
composition factors of the induced $L$-module $\frak h$, where $L$ is
the Levi component of some parabolic subgroup of $G$ containing a
fixed maximal torus $T$.

\proclaim{Corollary 8.5}
Assume that $\char k=p$ and $\mu _{\max} (\Omega _X)> 0$. There
exists a constant $B_G$ depending only on $G$ such that for every
principal $G$-bundle $E$ we have
$$\mu _{\max} (E(\frak g))\le \deg _{\HN} E +B_G\cdot {L_{\max}(\Omega _X)\over p}.$$
\endproclaim

\pr
Let $E_P$ be the Harder--Narasimhan filtration of $E$ and let
$L=P/R_u(P)$ be the Levi subgroup of $P$.  Since $E(\frak g)$ has a
filtration with quotients $E_L(V_S)$ for all possible shapes $S$, we have
$$\mu _{\max} (E(\frak g))\le \mathop{\max}_S{\mu _{\max} (E_L(V_S))}.\leqno (*)$$
Now let us recall that
$$\deg _{\HN}E=\deg E_P(\frak p)=\sum _{l(S)\ge 0}\deg E_L (V_S)$$
and $\deg E_L(V_S)$ are non-negative if $l(S)\ge 0$. Hence for any shape $S$ we have
$|\deg E_L(V_S)|\le \deg _{\HN}E$.
Since $E_L$ is semistable, Theorem 8.4 implies that there exists a constant $C_L(V_S)$
such that
$$\aligned
\mu _{\max}(E_L(V_S))-\mu _{\min}(E_L(V_S))\le &L _{\max}(E_L(V_S))-L
_{\min}(E_L(V_S)) \\&=L _{\max}(E_L(\gl V_S))\le C_L(V_S)\cdot {L_{\max}
(\Omega _X)\over p}.\\
\endaligned
$$ 
Since if we fix a maximal torus there are only finitely many
possible choices for $P$ and $L$, it follows that there exists $B_G$
such that $C_L(V_S)\le B_G$ for all possible $P$, $L$ and $S$. Then
$$\mu _{\max} (E_L(V_S))\le \mu (E_L(V_S))+B_G\cdot {L_{\max} (\Omega _X)\over
p}\le  {\deg _{\HN} E\over \dim V_S}+B_G\cdot {L_{\max} (\Omega _X)\over
p},$$
which by $(*)$ implies the required inequality, Q.E.D.

\medskip
{\sl Acknowledgements.}

This paper is a consequence of the author's stay at the Tata Institute
of Fundamental Research in Mumbai, India, in November 2003.  The
author would like to thank Professor V. Mehta for his kind invitation
and support. The author would like to thank people at the Tata
Institute, especially to Professors S. Subramanian and
A.~J. Parameswaran, for useful conversations. Finally, the author
would like to thank the referee for his remarks.

The author was partially supported by a Polish KBN grant (contract
number 1P03A03027) and by a subsidy of Professor A. Bialynicki-Birula
awarded by the Foundation for Polish Science.


\medskip
 
\Refs 
\widestnumber\key{SGA3} 
 
 
\ref\key AB
\by  M.~F.~Atiyah, R. Bott 
\paper The Yang--Mills equations over Riemann surfaces
\jour Phil. Trans. R. Soc. London A \vol 308 \yr 1982 \pages 523--615 
\endref

\ref\key ABS 
\by  H. Azad, M. Barry, G. Seitz 
\paper On the structure of parabolic subgroups 
\jour Comm. Algebra \vol 18 \yr 1990 \pages 551--562 
\endref 
 
\ref\key Be1
\by K.~A.~Behrend
\paper The Lefschetz trace formula for the moduli stack of principal bundles
\paperinfo  University of California, PhD thesis, 1991; available at 
{\tt http://www.math.ubc.ca/{$\sim$}behrend}
\endref


\ref\key Be2 
\by  K.~A.~Behrend 
\paper Semistability of reductive group schemes over curves 
\jour Math. Ann.\vol 301 \yr 1995 \pages 281--305 
\endref 
 
\ref\key BG 
\by  I. Biswas, T. Gomez 
\paper Restriction theorems for principal bundles 
\jour Math. Ann. \vol 327 \yr 2003 \pages 773--792 
\endref 

\ref\key SGA3
\by  M. Demazure, A. Grothendieck
\book Sch\'emas en groupes 
\bookinfo Lect. Notes in Math. \vol 151, 152, 153 \yr 1970 
\endref 

\ref\key HN 
\by  Y. Holla, M.~S.~Narasimhan 
\paper A generalisation of Nagata's theorem on ruled surfaces 
\jour Compositio Math.\vol 127 \yr 2001 \pages 321--332 
\endref 

\ref\key IMS
\by  S. Ilangovan, V.~B.~Mehta, A.~J.~Parameswaran 
\paper Semistability and semisimplicity in  representations of low height
in positive characteristic 
\inbook A tribute to C. S. Seshadri, Perspectives in Geometry and Representation Theory,
Hindustan Book Agency \yr 2003 \pages 271--282 
\endref 


\ref\key Ja
\by  J.~C.~Jantzen 
\book Representations of algebraic groups
\bookinfo Mathematical Surveys and Monographs {\bf 107}, 
Second edition
\yr 2003
\endref 
 
\ref\key KP
\by  H. Kraft, C. Procesi 
\book Classical invariant theory. A primer
\bookinfo preprint, 1996; available at {\tt http://www.math.unibas.ch}
\endref 

\ref\key La1 
\by  A. Langer 
\paper Semistable sheaves in positive characteristic 
\jour Ann. of Math. \vol  159\yr 2004 \pages 251--276
\endref 
 

\ref\key La2 
\by  A. Langer 
\paper Moduli spaces of sheaves in mixed characteristic 
\paperinfo Duke Math.~J. \vol 124 \yr 2004 \pages 571--586 
\endref 

\ref\key MS 
\by V.~B.~Mehta, S. Subramanian 
\paper On the Harder-Narasimhan filtration of principal bundles 
\inbook Algebra, arithmetic and geometry, Part I, II (Mumbai, 2000),   
Tata Inst. Fund. Res. Stud. Math. \bf{16}
\yr 2002 \pages 405--415 
\endref 

 
\ref\key RR 
\by S. Ramanan, A. Ramanathan 
\paper Some remarks on the instability flag 
\jour Tohoku Math.~J. \vol 36 \yr 1984 \pages 269--291 
\endref 
 

\ref\key Ra1 
\by A. Ramanathan 
\paper Deformations of principal bundles on the projective line 
\jour Invent. Math.\vol 71 \yr 1983 \pages 165--191 
\endref 

\ref\key Ra2 
\by A. Ramanathan 
\paper Moduli for principal bundles over algebraic curves: I 
\jour Proc. Indian Acad. Sci. (Math. Sci.)\vol 106 \yr 1996 \pages 301--328 
\endref 
 
\ref\key Ra3 
\by A. Ramanathan 
\paper Moduli for principal bundles over algebraic curves: II 
\jour Proc. Indian Acad.~Sci. (Math. Sci.)\vol 106 \yr 1996 \pages 421--449 
\endref 


\ref\key Sp
\by  T.~A.~Springer 
\book Linear algebraic groups
\bookinfo  Second Edition, Progress in Mathematics {\bf 9}, 
Birkh\"auser, 
\yr 1998
\endref 

\ref\key Su 
\by X. Sun 
\paper Remarks on semistability of $G$-bundles in positive characteristic 
\jour Compositio Math. \vol 119 \yr 1999 \pages 41--52 
\endref 
 
\endRefs
\enddocument